\documentclass[10pt,a4paper]{article}
\usepackage[margin=3.24cm]{geometry}
\usepackage{amsfonts,amsmath,bm,dsfont,mathrsfs}
\usepackage{amssymb,amsthm}
\usepackage{enumitem,cases}

\usepackage{tabularx} 
\usepackage{pdflscape}  
\usepackage{rotating}

\usepackage{graphicx}
\usepackage{subcaption}

\usepackage{epstopdf}
\usepackage{algorithm}
\usepackage{algpseudocode}
\usepackage{tikz}
\usepackage{booktabs}    
\usepackage{float,needspace}
\usepackage{multirow}

\usepackage[numbers]{natbib}
\usepackage{microtype}

\usepackage[dvipsnames]{xcolor}
\colorlet{xlinkcolor}{red!50!black}
\usepackage{hyperref}[6.83]
\definecolor{refblue}{RGB}{26,13,171}
\hypersetup{
  colorlinks = true,
  allcolors = refblue,
  urlcolor = BrickRed,
}

\newtheorem{assumption}{Assumption}

\newtheorem{proposition}{Proposition}
\newtheorem{corollary}{Corollary}

\numberwithin{equation}{section}
\numberwithin{table}{section}
\numberwithin{lemma}{section}
\numberwithin{example}{section}
\numberwithin{definition}{section}
\numberwithin{assumption}{section}
\numberwithin{theorem}{section}
\numberwithin{proposition}{section}
\numberwithin{corollary}{section}
\numberwithin{remark}{section}
\newcommand{\email}[1]{\href{mailto:#1}{\texttt{#1}}}

\newenvironment{MSCcodes}
 {\par\medskip\noindent\textbf{MSC2020:}\space\ignorespaces}
 {\par}

\setlist[description]{style=multiline,leftmargin=2em}
\setlist[itemize]{style=standard,leftmargin=2em}
\setlist[enumerate]{style=standard,leftmargin=2.2em,itemsep=3pt}

\makeatletter
\renewcommand{\maketitle}{%
 \begin{center}
 {\large\bfseries\@title\par}\vspace{0.7em}
 {\small\@author\par}\vspace{0.5em}
 {\footnotesize\authoraffiliations\par}\vspace{0.6em}
 {\small\@date\par}
 \end{center}\vspace{0.6em}}
\makeatother
\AddToHook{cmd/section/before}{\Needspace{9\baselineskip}}
\title{Recent Developments in HPR Methods for\\
Large-Scale Convex Composite Quadratic Programming}

\author{%
 Kaihuang Chen\textsuperscript{1}\quad
 Defeng Sun\textsuperscript{1}\quad
 Yancheng Yuan\textsuperscript{1}\quad
 Guojun Zhang\textsuperscript{1}\quad
 Xinyuan Zhao\textsuperscript{2}}
\newcommand{\authoraffiliations}{%
 \textsuperscript{1}Department of Applied Mathematics, The Hong Kong Polytechnic University,\par
 Hung Hom, Hong Kong\par
 \textsuperscript{2}Department of Mathematics, Beijing University of Technology, Beijing, P.R. China\par\vspace{0.3em}
 \email{kaihuang.chen@connect.polyu.hk}\quad \email{defeng.sun@polyu.edu.hk}\par
 \email{yancheng.yuan@polyu.edu.hk}\quad \email{guojun.zhang@connect.polyu.hk}\par
 \email{xyzhao@bjut.edu.cn}}
\hypersetup{pdfauthor={Kaihuang Chen, Defeng Sun, Yancheng Yuan, Guojun Zhang, Xinyuan Zhao}}
\date{1 November 2026}

\def \bR {\mathbb{R}}

\def\cK{{\cal K}}

\def\cM{{\cal M}}
\def\cC{{\cal C}}
\def\cT{{\cal T}}

\def\cH{{\cal H}}

\def\by{{\bar y}}
\def\bz{{\bar z}}
\def\bx{{\bar x}}
\def\bw{{\bar w}}

\def\bu{{\bar u}}

\newcommand{\cW}{\mathcal W}

\DeclareMathOperator{\Rank}{rank}

\DeclareMathOperator{\Diag}{Diag}

\begin{document}
\maketitle
\begin{abstract}
We present recent developments in Halpern Peaceman--Rachford (HPR)
methods for solving large-scale convex composite quadratic programming
problems. Within a unified HPR framework applied to the restricted
Wolfe dual, we develop several algorithmic variants based on different
choices of the semi-proximal operator $\mathcal T_1$ for the joint dual
subproblem. Specifically, an sGS-based construction leads to a
sequential $w$--$y$--$w$ update scheme, while a block-linearized
construction eliminates cross-block coupling and enables parallel
dual updates. To improve the latter, we further introduce a
coupling-aware construction that exploits a second spectral moment
to obtain a tighter row majorizer. We extend the framework to convex
quadratically constrained quadratic programming through factor-based
second-order cone reformulations. Additional computational enhancements
include equality Gram scaling and asynchronous CPU polishing.
We establish convergence results under fixed metrics and discuss
practical adaptive strategies.

\end{abstract}

\section*{Outline}
\begingroup
\renewcommand{\arraystretch}{1.3}
\begin{tabularx}{\textwidth}{@{}r@{\quad}Xr@{}}
\ref*{sec:Intro} & \hyperref[sec:Intro]{\nameref*{sec:Intro}} & \pageref*{sec:Intro} \\
\ref*{sec:framework} & \hyperref[sec:framework]{\nameref*{sec:framework}} & \pageref*{sec:framework} \\
\ref*{sec:sgs-qp} & \hyperref[sec:sgs-qp]{\nameref*{sec:sgs-qp}} & \pageref*{sec:sgs-qp} \\
\ref*{sec:parallel-hpr} & \hyperref[sec:parallel-hpr]{\nameref*{sec:parallel-hpr}} & \pageref*{sec:parallel-hpr} \\
\ref*{sec:qcqp-extension} & \hyperref[sec:qcqp-extension]{\nameref*{sec:qcqp-extension}} & \pageref*{sec:qcqp-extension} \\
\ref*{sec:restart} & \hyperref[sec:restart]{\nameref*{sec:restart}} & \pageref*{sec:restart} \\
\ref*{sec:scaling} & \hyperref[sec:scaling]{\nameref*{sec:scaling}} & \pageref*{sec:scaling} \\
\ref*{sec:polishing} & \hyperref[sec:polishing]{\nameref*{sec:polishing}} & \pageref*{sec:polishing} \\
\ref*{sec:conclusions} & \hyperref[sec:conclusions]{\nameref*{sec:conclusions}} & \pageref*{sec:conclusions} \\
\end{tabularx}
\endgroup
\clearpage

\section{Introduction}
\label{sec:Intro}
\href{https://github.com/PolyU-IOR/HPR-QP}{HPR-QP} applies the Halpern Peaceman--Rachford (HPR) method to large-scale convex quadratic optimization. Its basic iteration combines semi-proximal splitting of the restricted Wolfe dual with Halpern anchoring~\citep{parallel-hpr-framework}. This report brings together the algorithm and recent techniques for parallel dual updates, row conditioning and high-accuracy completion in a common notation. We emphasize the role of each component and cite the original sources for established analyses.

\subsection{Problem formulation}
We consider the convex composite quadratic program
\begin{equation}\label{QP-primal}
 \min_{x\in\mathbb R^n}
 \left\{\tfrac12\langle x,Qx\rangle+\langle c,x\rangle+\phi(x)
 \;\middle|\; Ax\in\mathcal K\right\},
\end{equation}
where $Q\succeq0$ is self-adjoint, $A:\mathbb R^n\to\mathbb R^m$ is linear, $\phi$ is proper, closed and convex, and $\mathcal K=[\ell,u]$ is a nonempty, possibly unbounded box. The main techniques concern the box-constrained QP
\begin{equation}\label{QP-primal-0}
 \min_{x\in\mathcal C}
 \left\{\tfrac12\langle x,Qx\rangle+\langle c,x\rangle
 \;\middle|\; Ax\in\mathcal K\right\},\qquad \mathcal C=[L,U],
\end{equation}
obtained with $\phi=\delta_{\mathcal C}$. Products by $Q$ need not use an explicitly formed matrix. The restricted Wolfe dual is
\begin{equation}\label{QP-dual}
 \begin{aligned}
 \min_{y,w,z}\quad &\tfrac12\langle w,Qw\rangle
       +\delta_{\mathcal K}^*(-y)+\phi^*(-z)\\
 \text{subject to}\quad &-Qw+A^*y+z=c,\qquad
 w\in\mathcal W:=\operatorname{Range}(Q).
 \end{aligned}
\end{equation}
The restriction makes $Q$ positive definite on $\mathcal W$. A full-space shadow vector represents $w$ without an explicit range projection~\citep{li2018qsdpnal,parallel-hpr-framework}.

\subsection{Notation and scope}
We write $A^*$ for the adjoint, $\|A\|$ for the spectral norm, and $\|v\|_M^2=\langle v,Mv\rangle$ for $M\succeq0$. The symbols $\delta_C$, $\Pi_C$ and $N_C$ denote the indicator, Euclidean projection and normal cone of a closed convex set $C$. For a proper closed convex function $f$, $f^*$ is its conjugate and $\operatorname{Prox}_f(v)=\arg\min_z\{f(z)+\tfrac12\|z-v\|^2\}$. Absent blocks are omitted; inverses act only on the indicated active spaces.

Bars denote proximal images and hats denote reflected iterates. The indices $k$, $r$ and $t$ denote total iterations, restart epochs and iterations within an epoch. The QP state is $u=(y,w,z,x)$, and $u_Q$ denotes its shadow version. The QCQP and polishing states are defined in their respective sections. Scaled data carry the subscript $\rm sc$; $D_x$ and $D_r$ denote variable and row scalings. For conic models written as $Ax-b\in\mathcal K$, the dual minimization includes $-\langle b,y\rangle$; thus $B_ix+h_i$ has offset $-h_i$ in that convention.

Section~\ref{sec:framework} presents the HPR method and its established convergence guarantees. Sections~\ref{sec:sgs-qp} and~\ref{sec:parallel-hpr} develop the sGS and block-linearized realizations, followed by the QCQP extension in Section~\ref{sec:qcqp-extension}. Sections~\ref{sec:restart}--\ref{sec:polishing} collect restart and penalty updates, scaling, and CPU polishing. Each technique is presented through its purpose, key formulas and guarantees.

The parallel method and public HPR-QP-C v0.1.1~\citep{hprqpc-release011} are distinct implementations. New metric derivations are included with their techniques; established proofs and public implementation details are available in the cited sources.

\Needspace{15\baselineskip}
\subsection{Software repositories}\label{intro:repositories}
The HPR project website is \url{https://hpr.top/}.
Table~\ref{intro:tab:repositories} lists the public GitHub repositories for
our Julia and C/CUDA implementations of HPR-QP.
Repository names are clickable. The linked default branches may differ
from the version-specific implementations discussed here.

\begin{table}[H]
\centering\small
\begin{tabularx}{\textwidth}{@{}ll>{\raggedright\arraybackslash}X@{}}
\toprule
Software & Execution & GitHub repository\\
\midrule
HPR-QP (Julia) & CPU, GPU &
\href{https://github.com/PolyU-IOR/HPR-QP}{\nolinkurl{PolyU-IOR/HPR-QP}}\\
HPR-QP-C (C/CUDA) & CPU, GPU &
\href{https://github.com/PolyU-IOR/HPR-QP-C}{\nolinkurl{PolyU-IOR/HPR-QP-C}}\\
\bottomrule
\end{tabularx}
\caption{Public HPR-QP implementations.}
\label{intro:tab:repositories}
\end{table}

HPR-QP-C also provides
\href{https://github.com/PolyU-IOR/HPR-QP-C/blob/main/python/README.md}{Python}
and \href{https://github.com/PolyU-IOR/HPR-QP-C/blob/main/julia/README.md}{Julia}
interfaces to the C/CUDA solver. The standalone Julia implementation in
Table~\ref{intro:tab:repositories} supports JuMP through MathOptInterface.

\section{The HPR method for QP}
\label{sec:framework}
\label{parallel:sec:model}

For $\sigma>0$, define the augmented Lagrangian of~\eqref{QP-dual} by
\begin{align}
 L_\sigma(y,w,z;x)
 ={}&\tfrac12\langle w,Qw\rangle+\delta_{\mathcal K}^*(-y)+\phi^*(-z)
       +\langle x,-Qw+A^*y+z-c\rangle\notag\\
 &+\tfrac\sigma2\|-Qw+A^*y+z-c\|^2.
 \label{parallel:eq:augmented-lagrangian}
\end{align}
In the QP case, $\phi=\delta_{\mathcal C}$. The following iteration applies to the composite model as well.

\subsection{Semi-proximal iteration}
\label{parallel:sec:framework}
Write $u=(y,w,z,x)$, $\bar u=(\bar y,\bar w,\bar z,\bar x)$ and $\mathcal U=\mathbb R^m\times\mathcal W\times\mathbb R^n\times\mathbb R^n$. Algorithm~\ref{alg:acc-pADMM-QP-0} is the common HPR method. Its defining choice is the self-adjoint positive semidefinite operator $\mathcal T_1$ in the joint $(y,w)$ subproblem.

\begin{algorithm}[H]
	\caption{An HPR method for solving the restricted Wolfe dual problem~\eqref{QP-dual}}
	\label{alg:acc-pADMM-QP-0}
	\begin{algorithmic}[1]
		\State \textbf{Input:} Choose a self-adjoint positive semidefinite linear operator \( \mathcal{T}_1 \) on \( \mathbb{R}^{m} \times \mathcal{W} \). Denote $u=(y,w,z,x)$ and $\bu=(\by,\bw,\bz,\bx)$. Let \( u^0 = (y^0, w^0, z^0, x^0) \in \mathcal{U} \), and set  \( \sigma > 0 \).
		\For{\( k = 0, 1, 2, \ldots \)}
		      \State {Step 1. $\displaystyle \bz^{k+1}=\underset{z \in \mathbb{R}^{n}}{\arg \min }\left\{L_\sigma(y^k, w^k,z; x^k)\right\} $;}
                \State {Step 2. $\displaystyle  \bx^{k+1}={x}^k+\sigma (-Qw^{k}+A^{*}y^{k}+\bz^{k+1}-c)$;}
                \State {Step 3. $\displaystyle (\by^{k+1},\bw^{k+1})=\underset{(y,w) \in \mathbb{R}^{m}\times \mathcal{W}}{\arg \min }\left\{\begin{aligned}&L_\sigma(y,w,\bz^{k+1};\bx^{k+1})\\[-2pt]&+\tfrac12\|(y,w)-(y^k,w^k)\|_{\mathcal T_1}^2\end{aligned}\right\}$;}
                \State {Step 4. $\displaystyle \hat{u}^{k+1}= 2 \bu^{k+1}-u^k$;}
	          \State {Step 5. $\displaystyle u^{k+1}=\frac{1}{k+2}u^{0}+\frac{k+1}{k+2}\hat{u}^{k+1}$;}
		\EndFor
	\end{algorithmic}
\end{algorithm}

The smooth Hessian of that subproblem is
\begin{equation}\label{parallel:eq:H-operator}
 \mathcal H=\begin{bmatrix}
  \sigma AA^*&-\sigma AQ\\
  -\sigma QA^*&Q+\sigma Q^2
 \end{bmatrix}.
\end{equation}
We require $\mathcal H+\mathcal T_1$ to be positive definite on $\mathbb R^m\times\mathcal W$. Sections~\ref{sec:sgs-qp} and~\ref{sec:parallel-hpr} give two ways to satisfy this condition and solve the same joint subproblem: an sGS sweep and a block-linearized update.

\subsection{Projection-based implementation}
The first two steps of Algorithm~\ref{alg:acc-pADMM-QP-0} reduce to one proximal evaluation. With
\begin{equation}\label{framework:eq:primal-projection}
 \begin{aligned}
 r_z^k&=x^k+\sigma(-Qw^k+A^*y^k-c),\\
 \bar x^{k+1}&=\operatorname{Prox}_{\sigma\phi}(r_z^k),\qquad
 \bar z^{k+1}=\sigma^{-1}(\bar x^{k+1}-r_z^k),\qquad
 \widehat x^{k+1}=2\bar x^{k+1}-x^k,
 \end{aligned}
\end{equation}
the box QP uses $\operatorname{Prox}_{\sigma\phi}=\Pi_{\mathcal C}$. The joint dual step remains determined by $\mathcal T_1$: Sections~\ref{sec:sgs-qp} and~\ref{sec:parallel-hpr} give its sequential and independent-block realizations. Their shadow implementations replace $w$ by $w_Q$ without changing $Qw$.

\subsection{Established convergence guarantees}

The fixed-metric results below are recalled from \cite[Propositions 2.1--2.2]{parallel-hpr-framework}. Their proofs follow the accelerated degenerate proximal-point analysis in \cite[Section 3.2 of the arXiv version]{sun2025accelerating}.

A primal-dual solution satisfies
\begin{equation}\label{eq:KKT-QP}
 Qw=Qx,\quad Ax\in\partial\delta_{\mathcal K}^*(-y),\quad
 x\in\partial\phi^*(-z),\quad -Qw+A^*y+z=c.
\end{equation}
Let $A_Q=[A^*\;-Q]$. Define the linear operator \( \mathcal{M} : \mathcal{U} \to \mathcal{U}\) as
\begin{equation}\label{def:M-QP}
\mathcal{M} =\begin{bmatrix}
\sigma A_Q^{*}A_Q
+ \mathcal{T}_1 & \quad  \mathbf{0}& \quad
A_{Q}^{*} \\
\mathbf{0}& \quad \mathbf{0} & \quad \mathbf{0}\\
A_{Q} & \quad \mathbf{0} &\quad  \frac{1}{\sigma} I_n
\end{bmatrix}.
\end{equation}
Here $I_n$ is the identity. The Schur complement of $\sigma^{-1}I_n$ is $\mathcal T_1\succeq0$, so $\mathcal M$ defines a seminorm, which may vanish on nonzero stored coordinates. We use the following assumptions.
\begin{assumption}\label{ass: CQ-QP}
There exists a vector $({y}^*,w^{*}, {z}^*,{x}^*)\in \mathcal{U}$ satisfying the KKT system \eqref{eq:KKT-QP}.
\end{assumption}
\begin{assumption}\label{ass: stronglyconvex}
The operator $\cT_1$ is a self-adjoint positive semidefinite linear operator such that $ \cT_1+\mathcal{H}$ is positive definite on $\mathbb{R}^{m}\times \mathcal{W}$.
\end{assumption}
Under Assumption~\ref{ass: stronglyconvex}, each subproblem in Algorithm~\ref{alg:acc-pADMM-QP-0} has a unique solution.
The KKT residual mapping is
\begin{equation}\label{def:KKT_residual}
\mathcal{R}(u) =
\begin{pmatrix}
Ax - \Pi_{\mathcal{K}}(Ax  - y) \\
Qw - Qx \\
x - \operatorname{Prox}_{\phi}(x - z) \\
c - A^{*}y + Qw - z
\end{pmatrix} \quad \forall\, u = (y, w, z, x) \in \mathbb{R}^{m} \times \mathcal{W} \times \mathbb{R}^{n} \times \mathbb{R}^{n}.
\end{equation}
\begin{proposition}[Convergence and last-iterate residual bounds]
\label{Prop:convergence}\label{prop:complexity-acc-pADMM}
Under Assumptions~\ref{ass: CQ-QP}--\ref{ass: stronglyconvex}, the proximal sequence $\{\bar u^{k+1}\}$ generated by Algorithm~\ref{alg:acc-pADMM-QP-0} converges to a primal-dual solution. For any KKT point $u^*$ and $R_0:=\|u^0-u^*\|_{\mathcal M}$,
\begin{align*}
 \|\bar u^{k+1}-u^k\|_{\mathcal M}&\leq\frac{R_0}{k+1},\\
 \|\mathcal R(\bar u^{k+1})\|&\leq
 \left(\frac{\sigma\|A_Q^*\|+1}{\sqrt\sigma}
       +\|\sqrt{\mathcal T_1}\|\right)\frac{R_0}{k+1},\qquad k\geq0.
\end{align*}
\end{proposition}
The corresponding dual objective-error bounds are given in
\cite[Proposition 2.2]{parallel-hpr-framework}.

The guarantees above assume a fixed penalty and metric. The practical restart and penalty rules below use them as motivation; unrestricted metric adaptation requires a separate convergence argument. The next two sections instantiate Algorithm~\ref{alg:acc-pADMM-QP-0} through different choices of $\mathcal T_1$.

\section{The sGS realization for QP}
\label{sec:sgs-qp}\label{sec:2}
The symmetric Gauss--Seidel (sGS) decomposition realizes the joint dual step of Algorithm~\ref{alg:acc-pADMM-QP-0} as a sequential $w$--$y$--$w$ sweep. This section gives the metric and one explicit implementation; proofs and extended baseline discussion are available in~\citep{parallel-hpr-framework,li2019block}.
\subsection{The sGS metric}
Choose $\mathcal S_y\succeq0$ and $\mathcal S_w\succeq0$ such that $\mathcal S_y+AA^*\succ0$. Define
\begin{equation}\label{def:S}
 \mathcal S=\sigma\Diag(\mathcal S_y,\mathcal S_w),
 \qquad
 \mathcal T_{1,\rm sGS}=\mathcal S+\widehat{\mathcal S}_{\rm sGS},
\end{equation}
where
\begin{equation}\label{def:sGS}
 \widehat{\mathcal S}_{\rm sGS}
 =\Diag\!\left(\sigma^2AQ(Q+\sigma Q^2+\sigma\mathcal S_w)^{-1}QA^*,\,0\right).
\end{equation}
The inverse acts on $\mathcal W$. The sGS decomposition theorem gives the following equivalence~\cite[Theorem 1]{li2019block}.
\begin{proposition}[sGS realization]\label{Prop:sGS}
For $\mathcal T_1=\mathcal T_{1,\rm sGS}$, Step 3 of Algorithm~\ref{alg:acc-pADMM-QP-0} is equivalent to
\begin{align*}
 \bar w^{k+1/2}&=\arg\min_{w\in\mathcal W}
 \{L_\sigma(y^k,w,\bar z^{k+1};\bar x^{k+1})+\tfrac\sigma2\|w-w^k\|_{\mathcal S_w}^2\},\\
 \bar y^{k+1}&=\arg\min_y
 \{L_\sigma(y,\bar w^{k+1/2},\bar z^{k+1};\bar x^{k+1})+\tfrac\sigma2\|y-y^k\|_{\mathcal S_y}^2\},\\
 \bar w^{k+1}&=\arg\min_{w\in\mathcal W}
 \{L_\sigma(\bar y^{k+1},w,\bar z^{k+1};\bar x^{k+1})+\tfrac\sigma2\|w-w^k\|_{\mathcal S_w}^2\}.
\end{align*}
Moreover, $\mathcal H+\mathcal T_{1,\rm sGS}\succ0$ on $\mathbb R^m\times\mathcal W$.
\end{proposition}
The proof specialized to HPR-QP is given in~\cite[Proposition 2.3]{parallel-hpr-framework}.

\subsection{Projection-based implementation}
For positive spectral upper bounds $\lambda_A\geq\lambda_{\max}(AA^*)$ and $\lambda_Q\geq\lambda_{\max}(Q)$, choose
\begin{align}
 \mathcal S_y&=\lambda_AI_m-AA^*,\label{def:S_y}\\
 \mathcal S_w&=Q(\lambda_QI_n-Q).\label{def:S_w}
\end{align}
Then
\begin{equation}\label{eq:T1}
 \mathcal T_{1,\rm sGS}
 =\Diag\!\left(\sigma\mathcal S_y+
 \frac{\sigma^2}{1+\sigma\lambda_Q}AQA^*,\,
 \sigma Q(\lambda_QI_n-Q)\right).
\end{equation}
Store $w_Q\in\mathbb R^n$ with $w=\Pi_{\mathcal W}w_Q$. Since $Qw=Qw_Q$ and $\langle w,Qw\rangle=\langle w_Q,Qw_Q\rangle$, the following formulas need only products by $Q$.

\begin{algorithm}[htbp]
\caption{sGS shadow sweep for Algorithm~\ref{alg:acc-pADMM-QP-0}}
\label{alg:acc-pADMM-QP-im}
\begin{algorithmic}[1]
\Require Current state $u_Q^k=(y^k,w_Q^k,z^k,x^k)$ and fixed $\sigma>0$
\State $r_z^k=x^k+\sigma(-Qw_Q^k+A^*y^k-c)$
\State $\bar x^{k+1}=\operatorname{Prox}_{\sigma\phi}(r_z^k)$,\quad $\widehat x^{k+1}=2\bar x^{k+1}-x^k$
\State $\displaystyle\bar w_Q^{k+1/2}=(\sigma\lambda_Qw_Q^k+\widehat x^{k+1})/(1+\sigma\lambda_Q)$
\State $r_y^k=A\bigl(\widehat x^{k+1}+\sigma(Qw_Q^k-Q\bar w_Q^{k+1/2})\bigr)-\sigma\lambda_Ay^k$
\State $\displaystyle\bar y^{k+1}=(\Pi_{\mathcal K}(r_y^k)-r_y^k)/(\sigma\lambda_A)$
\State $\displaystyle\bar w_Q^{k+1}=\bar w_Q^{k+1/2}+\frac{\sigma}{1+\sigma\lambda_Q}A^*(\bar y^{k+1}-y^k)$
\State Recover $\bar z^{k+1}=(\bar x^{k+1}-r_z^k)/\sigma$ when needed.
\State Apply reflection and Halpern averaging from Algorithm~\ref{alg:acc-pADMM-QP-0} to $u_Q$.
\end{algorithmic}
\end{algorithm}

For the box QP, $\operatorname{Prox}_{\sigma\phi}=\Pi_{\mathcal C}$. Matching the initial $(y,z,x)$ and setting $w^0=\Pi_{\mathcal W}w_Q^0$ gives the same $(\bar y,Q\bar w,\bar z,\bar x)$ as the restricted iteration. Convergence and the inherited $O(1/k)$ bounds follow from~\cite[Proposition 2.4 and Theorems 2.1--2.2]{parallel-hpr-framework}. If $Q=0$, omit both $w_Q$ updates; if $m=0$, omit the row stage. Restarted execution is described in Section~\ref{sec:restart}.

\subsection{Relationship with primal formulations}\label{sgs:sec:primal}
Two primal alternatives introduce either one constraint slack or separate constraint and objective variables:
\begin{align}
 \min_{x,s}\quad&\tfrac12\langle x,Qx\rangle+\langle c,x\rangle
       +\phi(x)+\delta_{\mathcal K}(s)
       &&\text{subject to }Ax=s,\label{QP-primal-4}\\
 \min_{x,s,v}\quad&\tfrac12\langle v,Qv\rangle+\langle c,x\rangle
       +\phi(x)+\delta_{\mathcal K}(s)
       &&\text{subject to }Ax=s,\ x=v.\label{QP-primal-5}
\end{align}
Their HPR iterations and comparison with the restricted-dual formulation are given in \cite[Appendix A and Section 4.2]{parallel-hpr-framework}. We retain these formulations for context and refer to that arXiv paper for the full algorithms.

\section{The block-linearized realization for QP}
\label{sec:parallel-hpr}

The second realization of Algorithm~\ref{alg:acc-pADMM-QP-0} uses block linearization. We choose $\mathcal T_1$ so that $\mathcal H+\mathcal T_1$ is block diagonal, making the $y$ and $w$ updates independent. The construction uses the normalized interaction between $A$ and $Q$ and a second spectral moment to control the added curvature.

Here parallel means that both blocks use the same old state and can be evaluated concurrently. The method does not require multiple GPUs. It is a separate branch from the public implementation used for Gram scaling and polishing; no combined configuration is assumed.

\subsection{A coupling-aware metric}
\label{parallel:sec:metric}

Choose spectral upper bounds
\begin{equation}
\label{parallel:eq:spectral-bounds}
\begin{aligned}
 &\lambda_A\geq\lambda_{\max}(AA^*),\qquad
 \lambda_Q\geq\lambda_{\max}(Q),\\
 & \lambda_{AQA^*}:=\lambda_{\max}(AQA^*),\qquad
 \lambda_{AQ^2A^*}:=\lambda_{\max}(AQ^2A^*),
\end{aligned}
\end{equation}
with $\lambda_A>0$ when $m>0$ and $\lambda_Q>0$ when $Q\neq0$.  Define
\begin{equation}
\label{parallel:eq:zeta-definition}
 \zeta_{AQ}
 :=\sqrt{\frac{\lambda_{AQA^*}}{\lambda_A\lambda_Q}}
\in[0,1].
\end{equation}
Since $Q^2\preceq\lambda_QQ$,
\begin{equation}
\label{parallel:eq:second-moment-order}
 0\leq\lambda_{AQ^2A^*}\leq\lambda_Q\lambda_{AQA^*}.
\end{equation}

\paragraph{Canceling the dual-block coupling.}

We construct $\cT_1$ directly from the coupled Hessian
\eqref{parallel:eq:H-operator}.  To obtain independent $y$- and $w$-updates, prescribe
the inexpensive block-diagonal target
\begin{equation}
\label{parallel:eq:diagonal-target}
 \cH+\cT_1=
 \begin{bmatrix}
  \sigma d_y\lambda_AI_m & 0\\
  0 & Q+\sigma d_w\lambda_QQ
 \end{bmatrix},
 \qquad d_y,d_w\geq1.
\end{equation}
The off-diagonal blocks of $\cT_1$ are then forced: they must be
$\sigma AQ$ and $\sigma QA^*$ to cancel the two coupling blocks of $\cH$.
Subtracting \eqref{parallel:eq:H-operator} from \eqref{parallel:eq:diagonal-target} gives the
two-scalar family
\begin{equation}
\label{parallel:eq:T1-dy-dw}
 \frac{1}{\sigma}\cT_1(d_y,d_w)
 =
 \begin{bmatrix}
  d_y\lambda_AI_m-AA^* & AQ\\
  QA^* & d_w\lambda_QQ-Q^2
 \end{bmatrix}.
\end{equation}

Define the $w$ coefficient by
\begin{equation}
\label{parallel:eq:dy-dw}
 \boxed{
 d_w(\tau):=
 \begin{cases}
  1,&\zeta_{AQ}=0,\\[1mm]
  \displaystyle 1+\frac{\zeta_{AQ}}{\tau},&\zeta_{AQ}>0,
 \end{cases}
 \qquad \tau>0.}
\end{equation}
For $\zeta_{AQ}>0$, this gives $d_w\lambda_Q>\lambda_Q$.  On $\cW$, the Schur-complement
condition for \eqref{parallel:eq:T1-dy-dw} is
\begin{equation}
\label{parallel:eq:T1-exact-schur}
 d_y\lambda_AI_m
 \succeq
 AA^*+AQ(d_w\lambda_QI_n-Q)^{-1}A^*.
\end{equation}

\paragraph{Exact curvature and a practical bound.}
For $\nu=d_w\lambda_Q>\lambda_Q$, the least scalar row curvature allowed by~\eqref{parallel:eq:T1-exact-schur} is
\begin{equation}\label{parallel:eq:exact-scalar-curvature}
 \ell^*(\nu)=\nu\,\lambda_{\max}\!\left(A(\nu I_n-Q)^{-1}A^*\right).
\end{equation}
Evaluating it requires shifted linear solves. The second-moment construction below replaces this resolvent by a bound based on fixed products with $A$, $Q$ and $A^*$.

For $d_w>1$,
\begin{equation}
\label{parallel:eq:two-moment-rational-bound}
 Q(d_w\lambda_QI_n-Q)^{-1}
 \preceq
 \frac{1}{d_w\lambda_Q}Q
 +\frac{1}{d_w\lambda_Q(d_w\lambda_Q-\lambda_Q)}Q^2.
\end{equation}
Consequently, the following coefficient satisfies
\eqref{parallel:eq:T1-exact-schur}:
\begin{equation}
\label{parallel:eq:dy-two-moment}
\boxed{
 d_y(\tau):=
 \begin{cases}
  1,
    &\zeta_{AQ}=0,\\[1mm]
  \displaystyle
  1+
  \frac{\lambda_{AQA^*}}
       {\lambda_A\lambda_Qd_w(\tau)}
  +\frac{\lambda_{AQ^2A^*}}
       {\lambda_A\lambda_Q^2d_w(\tau)(d_w(\tau)-1)},
    &\zeta_{AQ}>0.
 \end{cases}}
\end{equation}
For $\zeta_{AQ}>0$, substitution of \eqref{parallel:eq:dy-dw} gives
\begin{equation}
\label{parallel:eq:dy-two-moment-tau}
 d_y(\tau)
 =1+
 \frac{\zeta_{AQ}^2\tau}{\zeta_{AQ}+\tau}
 +\frac{\lambda_{AQ^2A^*}\tau^2}
 {\lambda_A\lambda_Q^2\zeta_{AQ}(\zeta_{AQ}+\tau)}.
\end{equation}
The one-moment coefficient $1+\zeta_{AQ}\tau$ is recovered in the
worst-aligned case
$\lambda_{AQ^2A^*}=\lambda_Q\lambda_{AQA^*}$.  More generally,
\begin{equation}
\label{parallel:eq:dy-improvement}
 (1+\zeta_{AQ}\tau)-d_y(\tau)
 =\frac{\tau^2
 (\lambda_Q\lambda_{AQA^*}-\lambda_{AQ^2A^*})}
 {\zeta_{AQ}(\zeta_{AQ}+\tau)\lambda_A\lambda_Q^2}
 \geq0.
\end{equation}
Thus, for the same $d_w$, the second-moment construction never decreases
the nominal $y$ step.

When $\zeta_{AQ}=0$, equivalently $AQ=0$, both coefficients equal one and
are independent of $\tau$.  With
\eqref{parallel:eq:dy-dw}--\eqref{parallel:eq:dy-two-moment}, choose
\begin{equation}
\label{parallel:eq:T1-parallel}
 \boxed{
 \cT_{1,\mathrm{par}}(\sigma,\tau)
 =\sigma
 \begin{bmatrix}
 d_y\lambda_A I_m-AA^* & AQ\\
 QA^* & d_w\lambda_QQ-Q^2
 \end{bmatrix}.}
\end{equation}

\begin{proposition}[Admissibility and separability]
\label{parallel:prop:metric-psd}
Suppose that \eqref{parallel:eq:spectral-bounds} are valid upper bounds.  Then, for
every $\sigma>0$ and $\tau>0$,
$\cT_{1,\mathrm{par}}(\sigma,\tau)\succeq0$.  Moreover,
\begin{equation}
\label{parallel:eq:diagonalized-H}
 \cH+\cT_{1,\mathrm{par}}
 =
 \begin{bmatrix}
   \sigma d_y\lambda_A I_m & 0\\
   0 & Q+\sigma d_w\lambda_QQ
 \end{bmatrix}.
\end{equation}
The right-hand side is positive definite on $\bR^m\times\cW$, and the
$y$- and $w$-subproblems in Algorithm~\ref{alg:acc-pADMM-QP-0} are independent.
\end{proposition}

\begin{proof}
If $\zeta_{AQ}=0$, then $AQ=0$ and $d_y=d_w=1$; hence
$\cT_{1,\mathrm{par}}/\sigma$ is the block-diagonal operator
$\Diag(\lambda_AI_m-AA^*,\lambda_QQ-Q^2)\succeq0$.

Suppose now that $\zeta_{AQ}>0$.  The lower-right block in
\eqref{parallel:eq:T1-parallel} is positive definite on $\cW$.  Its Schur complement
is positive semidefinite whenever \eqref{parallel:eq:T1-exact-schur} holds.
Inequality \eqref{parallel:eq:two-moment-rational-bound}, together with
$AA^*\preceq\lambda_AI_m$,
$AQA^*\preceq\lambda_{AQA^*}I_m$, and
$AQ^2A^*\preceq\lambda_{AQ^2A^*}I_m$, shows that
\eqref{parallel:eq:dy-two-moment} satisfies that condition.  Therefore
$\cT_{1,\mathrm{par}}\succeq0$.  Adding
\eqref{parallel:eq:T1-parallel} to \eqref{parallel:eq:H-operator} cancels the two cross blocks
and gives \eqref{parallel:eq:diagonalized-H}.
\end{proof}

\subsection{Projection-based implementation}
\label{parallel:sec:iteration}

Use the shadow state $u_Q=(y,w_Q,z,x)$ from Section~\ref{sec:sgs-qp}. Within epoch $r$, fix $\sigma_r,\tau_r>0$, evaluate $d_y,d_w$ at $\tau_r$, and set
\begin{equation}\label{parallel:eq:rho-gamma}
 \rho_y:=\sigma_r d_y\lambda_A,\qquad
 \gamma_w:=\sigma_r d_w\lambda_Q.
\end{equation}
Algorithm~\ref{parallel:alg:parallel-hpr} replaces the sequential sGS sweep. Its reflection, Halpern averaging and restart anchor follow the same convention as Section~\ref{sec:3}.

\begin{algorithm}[H]
\caption{Block-linearized shadow sweep for Algorithm~\ref{alg:acc-pADMM-QP-0}}
\label{parallel:alg:parallel-hpr}
\begin{algorithmic}[1]
\Require Current state $u_Q^{r,t}$ and fixed epoch parameters $\sigma_r,\tau_r$
\State $r_z^{r,t}=x^{r,t}+\sigma_r(-Qw_Q^{r,t}+A^*y^{r,t}-c)$
\State $\bar x^{r,t+1}=\Pi_{\mathcal C}(r_z^{r,t})$,\quad $\widehat x^{r,t+1}=2\bar x^{r,t+1}-x^{r,t}$
\State Compute the following two blocks independently:
\Statex $\displaystyle\qquad
 \bar w_Q^{r,t+1}=\frac{\gamma_w w_Q^{r,t}+\widehat x^{r,t+1}}{1+\gamma_w}$
\Statex $\displaystyle\qquad
 r_y^{r,t}=A\widehat x^{r,t+1}-\rho_y y^{r,t},\qquad
 \bar y^{r,t+1}=\frac{\Pi_{\mathcal K}(r_y^{r,t})-r_y^{r,t}}{\rho_y}$
\State Recover $\bar z^{r,t+1}=(\bar x^{r,t+1}-r_z^{r,t})/\sigma_r$ when needed.
\State Apply reflection and epoch-wise Halpern averaging to $u_Q$.
\end{algorithmic}
\end{algorithm}

The $y$ and $w_Q$ updates use the same reflected primal point and old dual state, so they can run concurrently. An ordinary iteration uses one product each by $Q$, $A$ and $A^*$, compared with two $Q$ products in the sGS path. The operator-dependency depth decreases from four layers to two. These work counts do not alone establish a shorter solve: setup, iteration count and termination accuracy also matter.

Table~\ref{parallel:tab:operator-work} summarizes the main operator work. The variable $z$ can be recovered when evaluating accuracy.

For completeness, the metric used by the fixed-point restart merit on the
full state $u_Q=(y,w_Q,z,x)$ is
\begin{equation}
\label{parallel:eq:M-parallel}
 \cM_{\mathrm{par}}(\sigma,\tau)=
 \begin{bmatrix}
  \sigma d_y\lambda_AI_m & 0 & 0 & A\\
  0 & \sigma d_w\lambda_QQ & 0 & -Q\\
  0 & 0 & 0 & 0\\
  A^* & -Q & 0 & \sigma^{-1}I_n
 \end{bmatrix}.
\end{equation}
For a displacement
$\delta u_Q=(\delta y,\delta w_Q,\delta z,\delta x)$,
\begin{align}
 \|\delta u_Q\|_{\cM_{\mathrm{par}}}^2
 ={}&\sigma d_y\lambda_A\|\delta y\|^2
 +\sigma d_w\lambda_Q
   \langle\delta w_Q,Q\delta w_Q\rangle \notag\\
 &+2\langle A^*\delta y-Q\delta w_Q,\delta x\rangle
 +\sigma^{-1}\|\delta x\|^2.
\label{parallel:eq:M-parallel-norm}
\end{align}

\subsection{Computational work}
The main computational advantage is a shorter dependency chain and one fewer product by $Q$ per ordinary iteration. Table~\ref{parallel:tab:operator-work} compares the two realizations for nonempty $A$ and general nonzero $Q$. Fewer products may reduce iteration cost, but total solve time also depends on setup, iteration count and accuracy.

\begin{table}[htbp]
\centering
\caption{Principal operator work per ordinary HPR-QP iteration.}
\label{parallel:tab:operator-work}
\begin{tabular}{@{}lcc@{}}
\toprule
Work item & sGS & Block-linearized \\
\midrule
Products by $Q$ & 2 & 1 \\
Products by $A$ and $A^*$ & 1 each & 1 each \\
Operator-dependency layers & 4 & 2 \\
\bottomrule
\end{tabular}
\end{table}

\section{Extension from QP to convex QCQP}
\label{sec:qcqp-extension}
The block-linearized QP method extends to convex quadratic constraints by adding factor-based cone blocks. We present the model, metric and projection updates together. Restart controls and scaling are collected in Sections~\ref{sec:restart} and~\ref{sec:scaling}, respectively.

\subsection{QCQP model and restricted Wolfe dual}
\label{parallel:sec:qcqp-extension}

Consider the convex QCQP
\begin{equation}
\label{parallel:eq:primal-qcqp}
 \begin{aligned}
  \min_{x\in\bR^n}\quad
    &\frac12\langle x,Qx\rangle+\langle c,x\rangle+\delta_{\cC}(x)\\
  \text{s.t.}\quad
    &Ax\in\cK,\\
    &\frac12\langle x,Q_i x\rangle+\langle c_i,x\rangle+d_i\leq0,
      \qquad i=1,\ldots,p,
 \end{aligned}
\end{equation}
where $Q_i\succeq0$.  Choose a thin factor $R_i$ such that
$Q_i=R_i^*R_i$, with $R_i:\bR^n\to\bR^{r_i}$, and let
\begin{equation*}
 \mathcal Q_i^{\rm r}
 :=\{(\xi,\eta,v)\in\bR\times\bR\times\bR^{r_i}
       \mid 2\xi\eta\geq\|v\|^2,\ \xi\geq0,\ \eta\geq0\}
\end{equation*}
be the rotated second-order cone.  Define
\begin{equation}
\label{parallel:eq:Bi-hi}
 B_i x:=
 \begin{bmatrix}
  0\\-c_i^*x\\R_i x
 \end{bmatrix},
 \qquad
 h_i:=
 \begin{bmatrix}
  1\\-d_i\\0
 \end{bmatrix}.
\end{equation}
Then
\begin{equation}
\label{parallel:eq:qc-rsoc-lift}
 B_i x+h_i=(1,-c_i^*x-d_i,R_i x)\in\mathcal Q_i^{\rm r}
\end{equation}
is equivalent to the $i$th quadratic inequality in
\eqref{parallel:eq:primal-qcqp}.  Stack
\begin{equation*}
 q:=\sum_{i=1}^p(r_i+2),
 \qquad
 B:=\begin{bmatrix}B_1\\ \vdots\\ B_p\end{bmatrix},
 \qquad
 h:=(h_1,\ldots,h_p),
 \qquad
 \mathcal Q^{\rm r}:=\mathcal Q_1^{\rm r}\times\cdots\times
                      \mathcal Q_p^{\rm r}.
\end{equation*}
Products by $B_i$ and $B_i^*$ use $c_i$ and $R_i$ directly, without forming $Q_i$. A minimal factor has $r_i=\Rank(Q_i)$, giving a cone of dimension $r_i+2$. A supplied factor with redundant rows remains valid, but its cone dimension is its actual row count plus two.

Here $s=(s_1,\ldots,s_p)$ is the cone multiplier; the symbol $s$ used locally in Section~\ref{sec:polishing} denotes a projected slack. The restricted Wolfe dual is
\begin{equation}
\label{parallel:eq:restricted-wolfe-qcqp}
 \begin{aligned}
  \min_{y,s,w,z}\quad
   &\frac12\langle w,Qw\rangle
    +\delta_{\cK}^*(-y)
    +\delta_{\mathcal Q^{\rm r}}^*(-s)
    +\langle h,s\rangle
    +\delta_{\cC}^*(-z)\\
  \text{s.t.}\quad
   &-Qw+A^*y+B^*s+z=c,
   \qquad w\in\cW.
 \end{aligned}
\end{equation}
Every rotated cone is self-dual, so
$\delta_{\mathcal Q^{\rm r}}^*(-s)$ is the indicator of
$s\in\mathcal Q^{\rm r}$.  As in Section~\ref{parallel:sec:iteration}, the implementation stores a
full-space shadow $w_Q\in\bR^n$ and uses only $Qw_Q$ and
$\langle w_Q,Qw_Q\rangle$; its projection onto $\cW$ represents the
restricted variable in \eqref{parallel:eq:restricted-wolfe-qcqp}.

\paragraph{HPR step for the QCQP dual.}
For $\sigma>0$, define
\begin{align}
\label{parallel:eq:augmented-lagrangian-qcqp}
 L_{\sigma}^{\rm qc}(y,s,w,z;x)
 :=\;&\frac12\langle w,Qw\rangle
 +\delta_{\cK}^*(-y)
 +\delta_{\mathcal Q^{\rm r}}^*(-s)
 +\langle h,s\rangle
 +\delta_{\cC}^*(-z) \notag\\
 &+\langle x,-Qw+A^*y+B^*s+z-c\rangle
 +\frac{\sigma}{2}\|-Qw+A^*y+B^*s+z-c\|^2.
\end{align}
Apply Algorithm~\ref{alg:acc-pADMM-QP-0} to $u_{\rm qc}=(y,s,w,z,x)$, replacing the joint $(y,w)$ minimization by a joint $(y,s,w)$ minimization of $L_\sigma^{\rm qc}$ with a positive semidefinite proximal operator on $\mathbb R^m\times\mathbb R^q\times\mathcal W$. The $z$ update, multiplier update, reflection and Halpern averaging retain their order. The explicit sweep is given below.

The Hessian of the smooth joint $(y,s,w)$ subproblem is
\begin{equation}
\label{parallel:eq:H-qcqp}
 \cH_{\rm qc}=\begin{bmatrix}
  \sigma AA^* & \sigma AB^* & -\sigma AQ\\
  \sigma BA^* & \sigma BB^* & -\sigma BQ\\
  -\sigma QA^* & -\sigma QB^* & Q+\sigma Q^2
 \end{bmatrix}.
\end{equation}

Before specifying the algorithm, choose a cone-compatible block-diagonal
majorizer
\begin{equation}
\label{parallel:eq:DB-abstract}
 D_B=\Diag(D_1,\ldots,D_p)\succeq BB^*,
 \qquad D_i\succ0,
\end{equation}
such that $D_i^{1/2}$ is an automorphism of $\mathcal Q_i^{\rm r}$.  The
simple baseline is
\begin{equation}
\label{parallel:eq:DB-baseline}
 D_B=\lambda_BI_q,
 \qquad
 \lambda_B\geq\lambda_{\max}(BB^*).
\end{equation}
It gives ordinary Euclidean rotated-cone projections.

Section~\ref{qcqp:sec:cone-majorizer} gives a tighter cone majorizer. The diagonal preprocessing in Section~\ref{sec:scaling} can be applied before its construction.

\subsection{The block-linearized metric}

Fix any majorizer satisfying \eqref{parallel:eq:DB-abstract}; the baseline choice is
\eqref{parallel:eq:DB-baseline}.  We majorize the three couplings $AB^*$, $AQ$, and
$BQ$ separately.  This pairwise construction is independent of a nesting
order and enlarges only the two blocks incident to each coupling.

Define
\begin{equation}
\label{parallel:eq:zeta-AB}
 \zeta_{AB}
 :=\left\|
 (\lambda_AI_m)^{-1/2}AB^*D_B^{-1/2}
 \right\|\in[0,1],
\end{equation}
and, when $Q\neq0$,
\begin{equation}
\label{parallel:eq:zeta-BQ}
 \zeta_{BQ}
 :=\sqrt{\frac{
 \lambda_{\max}\!\left(
 D_B^{-1/2}BQB^*D_B^{-1/2}
 \right)}{\lambda_Q}}
 \in[0,1].
\end{equation}
The factor $\zeta_{AQ}$ is the QP coupling in
\eqref{parallel:eq:zeta-definition}.  A factor is set to zero when either of its
blocks is absent.  Choose three independent shape parameters
$\tau_{AB},\tau_{AQ},\tau_{BQ}>0$.  Pairwise adaptation is described in Section~\ref{restart:sec:shapes}. Set
\begin{equation}
\label{parallel:eq:pairwise-alpha-qcqp}
 \begin{aligned}
  \alpha_y&:=\zeta_{AB}\tau_{AB}
             +d_y(\tau_{AQ})-1,\\
  \alpha_s&:=\frac{\zeta_{AB}}{\tau_{AB}}
             +\zeta_{BQ}\tau_{BQ},\\
  \alpha_w&:=d_w(\tau_{AQ})-1
             +\frac{\zeta_{BQ}}{\tau_{BQ}}.
 \end{aligned}
\end{equation}
Here $d_y(\tau_{AQ})$ and $d_w(\tau_{AQ})$ are exactly the QP coefficient
functions in \eqref{parallel:eq:dy-two-moment} and \eqref{parallel:eq:dy-dw}; hence their
$AQ$ contributions vanish when $AQ=0$.
The three diagonal coefficients are
\begin{equation}
\label{parallel:eq:dy-ds-dw-qcqp}
 d_y^{\rm qc}:=1+\alpha_y,
 \qquad
 d_s^{\rm qc}:=1+\alpha_s,
 \qquad
 d_w^{\rm qc}:=1+\alpha_w.
\end{equation}
Terms involving an absent block are omitted.  In particular, if $Q=0$,
then the $w_Q$ block and all $AQ$ and $BQ$ terms disappear.

The objective block remains $d_w^{\rm qc}\lambda_QQ$, exactly as in the QP method;
no eigenbasis, range projection, or per-coordinate objective weighting is
introduced.

\begin{proposition}[QCQP admissibility and separability]
\label{parallel:prop:qcqp-metric-psd}
Suppose that the spectral quantities in \eqref{parallel:eq:spectral-bounds},
\eqref{parallel:eq:zeta-AB}, and \eqref{parallel:eq:zeta-BQ} are valid upper bounds, and let
$D_B$ satisfy \eqref{parallel:eq:DB-abstract}.  Then, for every positive active
shape parameter,
\begin{equation}
\label{parallel:eq:T1-qcqp}
 \boxed{
 \cT_{1,\mathrm{par}}^{\rm qc}
 :=\sigma\left[
 \Diag\bigl(d_y^{\rm qc}\lambda_AI_m,\ d_s^{\rm qc}D_B,\ d_w^{\rm qc}\lambda_QQ\bigr)
 -
 \begin{bmatrix}A\\B\\-Q\end{bmatrix}
 \begin{bmatrix}A^*&B^*&-Q\end{bmatrix}
 \right]\succeq0.}
\end{equation}
Moreover,
\begin{equation}
\label{parallel:eq:HplusT-qcqp}
 \cH_{\rm qc}+\cT_{1,\mathrm{par}}^{\rm qc}
 =\Diag\bigl(
   \sigma d_y^{\rm qc}\lambda_AI_m,
   \sigma d_s^{\rm qc}D_B,
   (1+\sigma d_w^{\rm qc}\lambda_Q)Q
  \bigr).
\end{equation}
Consequently, $y$, all $s_i$, and $w_Q$ can be updated simultaneously.
\end{proposition}

\begin{proof}
Write the QCQP correction as the sum of the cone-block slack
\begin{equation*}
 \Diag\bigl(0,\ D_B-BB^*,\ 0\bigr)\succeq0
\end{equation*}
and the following three embedded two-block operators:
\begin{align*}
 E_{AB}&=
 \begin{bmatrix}
  \zeta_{AB}\tau_{AB}\lambda_AI_m&-AB^*\\
  -BA^*&\zeta_{AB}\tau_{AB}^{-1}D_B
 \end{bmatrix},\\
 E_{AQ}&=
 \begin{bmatrix}
  d_y(\tau_{AQ})\lambda_AI_m-AA^*&AQ\\
  QA^*&d_w(\tau_{AQ})\lambda_QQ-Q^2
 \end{bmatrix},\\
 E_{BQ}&=
 \begin{bmatrix}
  \zeta_{BQ}\tau_{BQ}D_B&BQ\\
  QB^*&\zeta_{BQ}\tau_{BQ}^{-1}\lambda_QQ
 \end{bmatrix}.
\end{align*}
The $AQ$ operator is positive semidefinite by
Proposition~\ref{parallel:prop:metric-psd}; the $AB$ and $BQ$ operators are positive
semidefinite by their defining normalized coupling bounds and the standard
two-block Schur-complement argument.  Together with the cone-block slack,
their sum is exactly the operator in \eqref{parallel:eq:T1-qcqp}, with the three
coefficients in \eqref{parallel:eq:pairwise-alpha-qcqp}; hence
\eqref{parallel:eq:T1-qcqp} is positive semidefinite.  Adding it to
\eqref{parallel:eq:H-qcqp} cancels all cross blocks and gives
\eqref{parallel:eq:HplusT-qcqp}.
\end{proof}

Let $F=[A^*\ B^*\ {-Q}]$ and
$D=\Diag(d_y^{\rm qc}\lambda_AI_m,d_s^{\rm qc}D_B,d_w^{\rm qc}\lambda_QQ)$.
The HPR metric, grouped as $((y,s,w),z,x)$, is
\begin{equation}\label{parallel:eq:M-qcqp}
 \mathcal M_{\rm qc}=\begin{bmatrix}
 \sigma D&0&F^*\\0&0&0\\F&0&\sigma^{-1}I_n
 \end{bmatrix}\succeq0.
\end{equation}
Its Schur complement is $\mathcal T_{1,\rm par}^{\rm qc}\succeq0$. Thus the fixed-metric guarantees of Section~\ref{sec:framework} apply to the enlarged model when a KKT point exists. In particular, the proximal sequence converges and its KKT residual has an $O(1/k)$ bound; this statement does not cover unrestricted adaptive changes.
For
$\delta u_{\rm qc}=(\delta y,\delta s,\delta w,\delta z,\delta x)$,
\begin{align}
\label{parallel:eq:M-qcqp-norm}
 \|\delta u_{\rm qc}\|_{\cM_{\rm qc}}^2
 ={}&\sigma d_y^{\rm qc}\lambda_A\|\delta y\|^2
 +\sigma d_s^{\rm qc}\langle\delta s,D_B\delta s\rangle
 +\sigma d_w^{\rm qc}\lambda_Q\langle\delta w,Q\delta w\rangle \notag\\
 &+2\langle A^*\delta y+B^*\delta s-Q\delta w,\delta x\rangle
 +\sigma^{-1}\|\delta x\|^2.
\end{align}
The implementation uses the same expression with the full-space shadow
$w_Q$, because the metric depends on that coordinate only through $Qw_Q$.

Equivalently, the exact primal-space admissibility condition is
\begin{equation}
\label{parallel:eq:qcqp-primal-condition}
 \boxed{
 \frac{1}{d_y^{\rm qc}\lambda_A}A^*A
 +\frac{1}{d_s^{\rm qc}}B^*D_B^{-1}B
 +\frac{1}{d_w^{\rm qc}\lambda_Q}Q
 \preceq I_n.}
\end{equation}
The final term is omitted when $Q=0$.  Condition
\eqref{parallel:eq:qcqp-primal-condition} is necessary and sufficient within the
prescribed block family, while the pairwise construction gives a scalable
sufficient way to choose its three block coefficients.

The pairwise construction is independent of a grouping order and enlarges only the two blocks incident to each coupling. In contrast, coupling a grouped block with a third block can enlarge both members even when only one interacts strongly.

\subsection{Block-diagonal cone majorizers}
\label{qcqp:sec:cone-majorizer}

Diagonal preprocessing alone does not supply $D_B\succeq BB^*$. The construction below combines scalings that preserve
each rotated cone with bounds on the interactions between constraint blocks.
For each retained quadratic constraint, choose
$\theta_i>0$ and define
\begin{equation}
\label{parallel:eq:Hi-theta}
 H_i(\theta_i):=\Diag(\theta_i^{-1},\theta_i,I_{r_i}).
\end{equation}
Both square roots of $H_i(\theta_i)$ are automorphisms of
$\mathcal Q_i^{\rm r}$.  In particular, the scaling used to construct the
majorizer acts on the complete affine block:
\begin{equation}
\label{parallel:eq:scaled-Bi-hi-majorizer}
 \widetilde B_i:=H_i(\theta_i)^{-1/2}B_i,
 \qquad
 \widetilde h_i:=H_i(\theta_i)^{-1/2}h_i,
 \qquad
 \widetilde B_ix+\widetilde h_i
 =H_i(\theta_i)^{-1/2}(B_ix+h_i).
\end{equation}
The weighted projection~\eqref{parallel:eq:s-update-qcqp} applies this map to the full affine residual, so both $B_i$ and $h_i$ are transformed without duplicate storage. For the canonical block, moreover,
\begin{equation}
\label{parallel:eq:Bi-Hi-Bi}
 B_i^*H_i(\theta_i)^{-1}B_i
 =Q_i+\frac{1}{\theta_i}c_ic_i^*.
\end{equation}
Thus the scaled metric still uses one ordinary rotated-cone projection.

Define
\begin{equation}
\label{parallel:eq:beta-i-theta}
 \beta_i:=\|\widetilde B_i\|^2.
\end{equation}
For the canonical unscaled block, \eqref{parallel:eq:Bi-Hi-Bi} gives
$\beta_i=\lambda_{\max}(Q_i+\theta_i^{-1}c_ic_i^*)$; after preprocessing,
the implementation evaluates \eqref{parallel:eq:beta-i-theta} from the scaled
$B_i$ data directly.
For the linear-memory construction, set
\begin{equation}
\label{parallel:eq:kappa-incidence}
 s_{iv}:=\|\widetilde B_ie_v\|,
 \qquad S_v:=\sum_hs_{hv},
 \qquad
 \kappa_i^{\rm inc}:=\beta_i+\sum_vs_{iv}(S_v-s_{iv}).
\end{equation}
If the sparse overlap graph is available, use instead
\begin{equation}
\label{parallel:eq:kappa-edge}
 \omega_{ij}:=\|\widetilde B_i\widetilde B_j^*\|,
 \qquad
 \kappa_i^{\rm edge}:=\beta_i+
     \sum_{j\in\mathcal N(i)}\omega_{ij}.
\end{equation}

\begin{proposition}[Scaled cone majorizer]
\label{parallel:prop:DB-majorizes-BBt}
Let $\kappa_i$ be either coefficient in
\eqref{parallel:eq:kappa-incidence}--\eqref{parallel:eq:kappa-edge}, with every overlap edge
included for the edge construction.  Then
\begin{equation}
\label{parallel:eq:DB-majorizer}
 D_B:=\Diag\bigl(\kappa_1H_1(\theta_1),\ldots,
                  \kappa_pH_p(\theta_p)\bigr)
 \succeq BB^*.
\end{equation}
\end{proposition}

\begin{proof}
After scaling by $\Diag(H_1^{-1/2},\ldots,H_p^{-1/2})$, block $i$ has
diagonal norm $\beta_i$, while each cross-block norm is bounded by
$\omega_{ij}$ or $\sum_vs_{iv}s_{jv}$.  Block Gershgorin and
$2ab\leq a^2+b^2$ give the stated diagonal margins.
\end{proof}

The weighted cone update requires $\kappa_i>0$. A zero operator block must be handled separately or assigned a positive valid majorizer.

The offset $h_i$ is absent from \eqref{parallel:eq:beta-i-theta}--
\eqref{parallel:eq:DB-majorizer} because it does not enter the Hessian $BB^*$; it is
nevertheless transformed in every cone evaluation according to
\eqref{parallel:eq:qcqp-scaled-affine-cone} and
\eqref{parallel:eq:scaled-Bi-hi-majorizer}.

Set $D_i=\kappa_iH_i(\theta_i)$ and recompute $\zeta_{AB},\zeta_{BQ}$. The cone update is unchanged. Substitution into~\eqref{parallel:eq:qcqp-primal-condition} gives
\begin{equation}\label{parallel:eq:qcqp-primal-condition-scaled}
 \frac{A^*A}{d_y^{\rm qc}\lambda_A}
 +\frac1{d_s^{\rm qc}}\sum_i\frac{B_i^*H_i(\theta_i)^{-1}B_i}{\kappa_i}
 +\frac{Q}{d_w^{\rm qc}\lambda_Q}\preceq I_n.
\end{equation}
This operator form remains valid after affine preprocessing.

\subsection{Projection-based implementation}

For the state $u_{Q,\mathrm{qc}}=(y,s,w_Q,z,x)$, first form
\begin{equation}
\label{parallel:eq:rz-qcqp}
 r_z^k=x^k+\sigma
 \left(-Qw_Q^k+A^*y^k+B^*s^k-c\right),
\end{equation}
followed by
\begin{equation}
\label{parallel:eq:xhat-qcqp}
 \bar x^{k+1}=\Pi_{\cC}(r_z^k),
 \qquad
 \widehat x^{k+1}=2\bar x^{k+1}-x^k,
 \qquad
 \bar z^{k+1}=\frac{\bar x^{k+1}-r_z^k}{\sigma}.
\end{equation}
The following blocks then depend only on $\widehat x^{k+1}$ and their old
shadow values and can therefore be evaluated concurrently.

For the linear block, let $\rho_y:=\sigma d_y^{\rm qc}\lambda_A$ and set
\begin{equation}
\label{parallel:eq:y-update-qcqp}
 r_y^k=A\widehat x^{k+1}-\rho_y y^k,
 \qquad
 \bar y^{k+1}
 =\frac{\Pi_{\cK}(r_y^k)-r_y^k}{\rho_y}.
\end{equation}

For cone $i$, set
\begin{equation}\label{parallel:eq:s-update-qcqp}
\begin{aligned}
 r_{s_i}^k&=B_i\widehat x^{k+1}+h_i-\sigma d_s^{\rm qc}D_is_i^k,\\
 \bar s_i^{k+1}
 &=\frac1\sigma(d_s^{\rm qc}D_i)^{-1/2}
 \left[\Pi_{\mathcal Q_i^{\rm r}}\bigl((d_s^{\rm qc}D_i)^{-1/2}r_{s_i}^k\bigr)
 -(d_s^{\rm qc}D_i)^{-1/2}r_{s_i}^k\right].
\end{aligned}
\end{equation}
Thus every weighted cone resolvent uses one ordinary rotated-SOC
projection.

For the objective block, let $\gamma_w:=\sigma d_w^{\rm qc}\lambda_Q$ and retain
the QP update
\begin{equation}
\label{parallel:eq:w-update-qcqp}
 \boxed{
 \bar w_Q^{k+1}
 =\frac{\gamma_w w_Q^k+\widehat x^{k+1}}{1+\gamma_w}.}
\end{equation}
No projection onto $\cW$ and no linear solve are required.  Apply the common reflection and Halpern averaging; restart controls are given in Section~\ref{sec:restart}.

\begin{corollary}[Exact reduction to QP]\label{parallel:cor:qcqp-reduces-qp}
If $p=0$, with matched parameter sequences and restart decisions, the QCQP iteration reduces exactly to Algorithm~\ref{parallel:alg:parallel-hpr}.
\end{corollary}
\begin{proof}
The cone blocks disappear, $\zeta_{AB}=\zeta_{BQ}=0$, and identifying $\tau_{AQ}=\tau$ gives $d_y^{\rm qc}=d_y(\tau)$ and $d_w^{\rm qc}=d_w(\tau)$. The remaining updates, metric and stopping test are those of the QP method.
\end{proof}

\subsection{Connection with standard second-order cones}

\label{sec:soc-support}

The same projection-based viewpoint applies to standard second-order-cone domains. This gives a common representation for conic QPs and factor-based QCQPs, as used in HPR-QP-C~\citep{hprqpc-release011}.

\paragraph{Cone domains and projections.}

For $d\geq1$, let
$\mathcal Q^d=\{(t,v):t\geq\|v\|_2\}$ be the Lorentz cone, with its scalar
coordinate first. The native model is
\begin{equation}\label{soc:eq:model}
 \min_{x\in\mathcal X}\;\tfrac12 \langle x,Qx\rangle+\langle c,x\rangle+c_0,
 \qquad Ax-b\in\mathcal K,
\end{equation}
where $\mathcal X$ combines box variables and variable-cone blocks, and
$\mathcal K$ combines linear intervals and affine row-cone blocks.
For $r=\|v\|_2$,
the block projection is
\begin{equation}\label{soc:eq:projection}
 \Pi_{\mathcal Q^d}(t,v)=
 \begin{cases}
 (t,v),&r\leq t,\\
 0,&r\leq-t,\\
 \tfrac12(t+r)(1,v/r),&\text{otherwise}.
 \end{cases}
\end{equation}
For $d=1$ this is scalar nonnegative projection. Product-domain projections
replace the interval/box projections in the ordinary two-$w$ HPR map;
the reflection and Halpern averaging are retained.

Affine offsets are transformed together with their row maps, as described in Section~\ref{sec:scaling}.

\paragraph{Sparse-factor reformulation.}

Suppose a convex quadratic row is supplied as
\begin{equation}\label{soc:eq:quadratic}
 g_i(x)=\tfrac12\|R_i x\|_2^2+\langle c_i,x\rangle+d_i\leq0,
\end{equation}
with a sparse rectangular factor $R_i$ having $r_i$ rows. Set
$t_i=-\langle c_i,x\rangle-d_i$. Then
\begin{equation}\label{soc:eq:lift}
 g_i(x)\leq0\quad\Longleftrightarrow\quad
 \begin{pmatrix}t_i+1\\t_i-1\\\sqrt2R_i x\end{pmatrix}
 \in\mathcal Q^{r_i+2}.
\end{equation}
Indeed, the difference between the squared leading coordinate and squared
tail norm is $4t_i-2\|R_i x\|_2^2$. Cone membership implies
$t_i\geq\tfrac12\|R_i x\|_2^2$; conversely that inequality gives both the
nonnegative leading coordinate and the cone inequality. This proves the
equivalence, including zero-rank factors.

The reformulation adds affine cone rows without adding optimization variables. It uses the factor $R_i$ directly and avoids forming $R_i^*R_i$; neither diagonality nor full row rank is required. The feasible set is preserved exactly, while finite-precision accuracy is assessed on the conic model.

\section{Restart and penalty updates}
\label{sec:restart}
Restarting and parameter adaptation aim to improve practical progress while retaining the same HPR sweep within each epoch. Their role is to choose a useful anchor and balance the primal and dual components. Baseline restart and penalty strategies are described in~\citep{parallel-hpr-framework}.

\subsection{Restarting the Halpern iteration}\label{sec:3}
Within an epoch, the Halpern anchor and metric remain fixed. Restarting replaces the anchor by an accepted proximal point and resets the averaging weights. Progress can be monitored through the metric displacement between the current iterate and its proximal image. A restart is motivated by sufficient progress, stagnation after improvement, or an overly long epoch~\cite[Sections 3.2--3.3]{parallel-hpr-framework}.

\subsection{Penalty and metric balance}\label{restart:sec:shapes}
The penalty $\sigma$ controls the relative weighting of primal and dual movement. The block-linearized method also has shape parameters that distribute curvature between coupled blocks. These parameters can be estimated from the displacement since the current anchor and updated at a restart. In QCQP, the same principle applies to the three pairwise couplings. After an update, progress measures must use the new metric consistently.


\section{Scaling techniques}
\label{sec:scaling}
Scaling prepares the row and column geometry before iteration. Diagonal equilibration preserves the variable and cone domains; equality Gram scaling further changes the row metric through a reusable linear solve. Both use the same coordinate convention and recover stopping quantities in the original model.

\subsection{Coordinate transformations}
Choose positive diagonal variable and ordinary-row factors $D_x,D_r$.
With $x=D_xx_{\rm sc}$, the ordinary data transform as
\begin{equation}\label{parallel:eq:qcqp-scaled-linear-data}
 \begin{aligned}
 Q_{\rm sc}&=D_xQD_x,&c_{\rm sc}&=D_xc,&A_{\rm sc}&=D_rAD_x,\\
 \mathcal C_{\rm sc}&=D_x^{-1}\mathcal C,
 &\mathcal K_{\rm sc}&=D_r\mathcal K.
 \end{aligned}
\end{equation}

For box QPs, these transformations preserve the model after recovery by $x=D_xx_{\rm sc}$. For QCQP and native conic models, the cone blocks require the additional restrictions below. The Gram option belongs to public HPR-QP-C v0.1.1.

\subsection{Cone-preserving diagonal equilibration}
Scaling transforms the complete affine cone map; metric construction majorizes only its linear part. Thus $h_i$ participates in equilibration and projections, but not in the Hessian $BB^*$.

\label{parallel:subsub:qcqp-scaling}\label{parallel:subsec:qcqp-enhancements}
For each rotated cone, we adapt Curtis--Reid scaling~\citep{scaling-curtis1972} using the automorphism
\begin{equation}\label{parallel:eq:qcqp-cr-automorphism}
 E_i(a_i,b_i)=\Diag(e^{a_i-b_i},e^{a_i+b_i},e^{a_i}I_{r_i}).
\end{equation}
It scales both $2\xi\eta$ and $\|v\|^2$ by $e^{2a_i}$, preserving cone membership. Subsequent Ruiz~\citep{ruiz2001scaling} and Pock--Chambolle~\citep{pock2011diagonal} passes use a common factor $\gamma_i>0$ within each cone:
\begin{equation}\label{parallel:eq:Gamma-cone-scaling}
 \Gamma=\Diag(\gamma_1I_{r_1+2},\ldots,\gamma_pI_{r_p+2}).
\end{equation}
After the Curtis--Reid transformation, each scalar cone pass applies
\begin{equation}\label{parallel:eq:qcqp-scaled-affine-cone}
 [B_{i,\rm sc}\ \ h_{i,\rm sc}]=\gamma_i[B_iD_x\ \ h_i].
\end{equation}
Scaling the full affine block gives the exact equivalence
\begin{equation}\label{parallel:eq:qcqp-scaling-equivalence}
 B_ix+h_i\in\mathcal Q_i^{\rm r}
 \quad\Longleftrightarrow\quad
 B_{i,\rm sc}x_{\rm sc}+h_{i,\rm sc}\in\mathcal Q_i^{\rm r}.
\end{equation}

The scaling sequence applies joint Curtis--Reid equilibration, followed by block Ruiz and Pock--Chambolle scaling. All stages respect the cone geometry and include the affine offsets in cone-row scaling. The variable factors are shared across the objective and constraint operators. When $p=0$, the cone transformations disappear and the procedure reduces to QP scaling. Solutions and stopping quantities are recovered in the original coordinates.

\paragraph{Standard Lorentz blocks.}

Each Ruiz/Pock scaling pass uses a common positive factor within a Lorentz cone to preserve its geometry. Row factors also scale affine offsets.
The same restriction applies to cone-variable blocks.

For an affine row map $Ax-b$, also set $b_{\rm sc}=D_rb$. Dual variables satisfy
$y_{\rm sc}=D_r^{-*}y$, preserving
$A_{\rm sc}^*y_{\rm sc}=D_xA^* y$. Cones remain standard when the
corresponding blocks of $D_x,D_r$ are positive scalars. Results and stopping
residuals are recovered in original native units. Gram whitening in
Section~\ref{sec:gram-scaling} mixes equality rows only and respects this
restriction on the other cone rows.

\subsection{Equality Gram scaling}
\label{sec:gram-scaling}
Equality Gram scaling changes the geometry of the row update by whitening the equality block. The resulting step uses a reusable Gram solve and a factor-two majorizer for its coupling to the remaining rows. This section describes the formulation in the HPR-QP-C implementation~\citep{hprqpc-release011,gram-coordinate-note}.

\paragraph{Whitening the equality rows.}
Partition the constraint operator into equality and remaining rows:
\begin{equation}\label{gram:eq:model}
 \min_{x\in\mathcal C}\;\tfrac12\langle x,Qx\rangle+\langle c,x\rangle+c_0,
 \qquad A_Ex=b_E,\qquad A_Rx-b_R\in\mathcal K_R,
 \qquad A=\begin{pmatrix}A_E\\A_R\end{pmatrix}.
\end{equation}
Here $\mathcal K_R$ is a product of interval and Lorentz blocks. Assume $A_E$ has full row rank, and set
\begin{equation}\label{gram:eq:whitening}
 \begin{aligned}
 G&=A_EA_E^*=L_GL_G^*\succ0,\qquad
 \beta>0,\quad\beta\geq\|A_R\|^2,\\
 R_G&=\operatorname{diag}(L_G^{-1},\beta^{-1/2}I).
 \end{aligned}
\end{equation}
The factor $L_G$ may be a Cholesky factor. The row bound $\beta$ is distinct from the HPR penalty $\sigma$; the result below requires a valid upper bound on $\|A_R\|^2$.

With $\widetilde A=R_GA$ and $\widetilde y=R_G^{-*}y$,
\begin{equation}\label{gram:eq:pairing}
 \widetilde A^*\widetilde y=A^*y,
 \qquad \langle\widetilde y,R_Gv\rangle=\langle y,v\rangle.
\end{equation}
Right-hand sides, affine offsets and finite interval endpoints are transformed with their rows. Invertible equality-row mixing preserves the equality constraints, and a common positive scalar on each Lorentz block preserves its cone. The primal variable, objective and box are unchanged. The implementation applies $G^{-1}$ without forming a dense whitening matrix or inverse.

\paragraph{The common majorizer.}
\begin{proposition}\label{gram:prop:majorizer}
Under~\eqref{gram:eq:whitening}, the row majorizer
\begin{equation}\label{gram:eq:majorizer}
 M_A=\operatorname{diag}(2G,2\beta I)\succeq AA^*,
 \qquad \|R_GA\|^2\leq2.
\end{equation}
\end{proposition}
\begin{proof}
The decomposition
\begin{equation}\label{gram:eq:decomposition}
 M_A-AA^*=
 \begin{pmatrix}A_E\\-A_R\end{pmatrix}
 \begin{pmatrix}A_E^*&-A_R^*\end{pmatrix}
 +\operatorname{diag}\bigl(0,2(\beta I-A_RA_R^*)\bigr)\succeq0
\end{equation}
proves the first claim. Multiplication by $R_G$ and $R_G^*$ gives $2I\succeq R_GAA^*R_G^*$, proving the second.
\end{proof}
The factor two controls the coupling between equality and remaining rows. With $A_E=A_R=[1]$ and $\beta=1$, the matrix $\operatorname{diag}(G,\beta I)-AA^*$ has a negative eigenvalue, so dropping the factor is invalid in general. With equality rows alone, $M_A=G$ would suffice mathematically; the documented implementation retains the common factor-two rule.

Ordinary row linearization uses $\sigma(\lambda_AI-AA^*)$. The Gram choice replaces this by $\sigma(M_A-AA^*)$. Equivalently, it is scalar-majorized HPR with bound two in whitened coordinates. The row majorizer $M_A$ is distinct from the full-state HPR metric $\mathcal M$.

\subsection{The Gram row update and controller geometry}
Let $v$ be the primal argument at the row stage of the sGS sweep. In Algorithm~\ref{alg:acc-pADMM-QP-im}, it is $v=\widehat x+\sigma(Qw_Q-Q\bar w_Q^{1/2})$, with iteration indices suppressed. Transforming the row step back to the original coordinates gives
\begin{equation}\label{gram:eq:update}
 \begin{aligned}
 Gq&=A_Ev-b_E,&\qquad \bar y_E&=y_E-\frac{q}{2\sigma},\\
 \rho_R&=2\sigma\beta,&
 \bar y_R&=\frac{\Pi_{\mathcal K_R}(A_Rv-b_R-\rho_Ry_R)
                     -(A_Rv-b_R-\rho_Ry_R)}{\rho_R}.
 \end{aligned}
\end{equation}
The equality step uses a Gram solve; the remaining rows retain a projection with the corresponding scalar $\rho_R$. The objective-block sweep, reflection and Halpern averaging retain their HPR form.

With the other dual block fixed, an unregularized equality ADMM minimization would give
\[
 G(\bar y_E-y_E)=-\frac{A_Ev-b_E}{\sigma}.
\]
Equation~\eqref{gram:eq:update} takes half this correction: the equality augmented-Lagrangian curvature $\sigma G$ is supplemented by the same proximal curvature. Removing that factor requires a different justified majorizer or block scheme.

\paragraph{Composition with diagonal scaling.}
Use the consistent variable and row scaling convention
\[
 x=D_xx_{\rm sc},\qquad A_{\rm sc}=D_rAD_x.
\]
The factors are positive diagonal and constant within each cone block. Let $D_{r,E}$ be the equality-row part of $D_r$, and define
\[
 A_{E,\rm sc}=D_{r,E}A_ED_x,\qquad G_{\rm sc}=A_{E,\rm sc}A_{E,\rm sc}^*.
\]
Construct $R_{G,\rm sc}$ from this scaled model. The combined operator is $R_{G,\rm sc}D_rAD_x$. Column scaling changes $G_{\rm sc}$, and remaining-row scaling changes the bound $\beta_{\rm sc}$.

For fixed $D_x$, exact whitening cancels invertible equality-row scaling in the primal-space projector. With $F_E=A_ED_x$,
\begin{equation}\label{gram:eq:row-invariance}
 A_{E,\rm sc}^*(A_{E,\rm sc}A_{E,\rm sc}^*)^{-1}A_{E,\rm sc}
 =F_E^*(F_EF_E^*)^{-1}F_E.
\end{equation}
This follows by substituting $(D_{r,E}F_EF_E^*D_{r,E}^*)^{-1}=D_{r,E}^{-*}(F_EF_E^*)^{-1}D_{r,E}^{-1}$. This identity describes the interaction of the two scalings in exact arithmetic.

Restart and penalty objectives use the same row geometry:
\begin{equation}\label{gram:eq:control}
 \langle\Delta y,M_A\Delta y\rangle
 =2\|A_E^*\Delta y_E\|^2+2\beta\|\Delta y_R\|^2
 =2\|\Delta\widetilde y\|^2.
\end{equation}
Thus restart and penalty measures should use the same row geometry as the iteration. Stopping tests remain in original units. Coordinate equivalence preserves the fixed-metric analysis under its assumptions; adaptive metric changes require a separate argument.

\subsection{Computational role}
The fixed Gram matrix can be factored once and reused across iterations and penalty changes. This trades additional setup, storage and linear solves for better row conditioning. Its use requires a full-row-rank equality block and an accurate inverse action. When these conditions or the computational tradeoff are unfavorable, ordinary row linearization remains available. Implementation details are given in~\citep{gram-coordinate-note,hprqpc-release011}.

\section{Asynchronous CPU polishing}
\label{sec:polishing}
Polishing targets high-accuracy completion from a promising HPR point. A frozen generalized-Jacobian correction is computed on the CPU from a copied primal--dual point while GPU HPR continues~\citep{hprqpc-release011}.

\subsection{A frozen correction equation}

Write the unscaled native problem as
\begin{equation}\label{polish:eq:model}
 \min_{x\in\mathcal C}\;
 f(x)=\tfrac12 \langle x,Qx\rangle+\langle c,x\rangle+c_0,
 \qquad Ax-b\in\mathcal K,
\end{equation}
where $Q\succeq0$, $\mathcal C$ is a box, and $\mathcal K$ is a product of
intervals and standard Lorentz cones. For QCQP input, ``native'' refers to
the lifted conic model in Section~\ref{sec:soc-support}.
Use the solver's signs
\[
 Qx+c-A^* y-z=0,\qquad
 -z\in N_{\mathcal C}(x),\qquad -y\in N_{\mathcal K}(Ax-b),
\]
where $N$ is the convex normal cone. For the polishing state $u_{\rm pol}=(x,y)$, define
\begin{equation}\label{polish:eq:residual}
 \begin{aligned}
 g&=Qx+c-A^* y, &p&=\Pi_{\mathcal C}(x-g),\\
 s&=\Pi_{\mathcal K}(Ax-b-y), &
 F(u_{\rm pol})&=\begin{pmatrix}r\\t\end{pmatrix}
       =\begin{pmatrix}x-p\\Ax-b-s\end{pmatrix}.
 \end{aligned}
\end{equation}
The projection
characterization $v-\Pi_C(v)\in N_C(\Pi_C(v))$ shows that $F(u_{\rm pol})=0$ implies
the native KKT conditions with $z=g$; conversely, a KKT point has $F(u_{\rm pol})=0$.

At the primal--dual pair $u_{\rm pol}^0$ copied from HPR, select and freeze generalized derivatives
$E$ and $S$ of the box and row projections in~\eqref{polish:eq:residual}.
The common frozen action is
\begin{equation}\label{polish:eq:jacobian}
 J_0\begin{pmatrix}\Delta x\\\Delta y\end{pmatrix}
 =\begin{pmatrix}
 (I-E+EQ)\Delta x-EA^*\Delta y\\
 (I-S)A\Delta x+S\Delta y
 \end{pmatrix}.
\end{equation}
Each correction approximately solves
\begin{equation}\label{polish:eq:chord}
 J_0d_k=-F(u_{\rm pol}^k),\qquad u_{\rm pol}^{k+1}=u_{\rm pol}^k+d_k.
\end{equation}
The projections and $F(u_{\rm pol}^k)$ are reevaluated, but $J_0$ is never refreshed
within an episode. This frozen generalized-Jacobian (chord) correction is
a local completion step.

\subsection{Candidate construction and accuracy}

The candidate associated with a polishing state is
\begin{equation}\label{polish:eq:candidate}
 x_c=p,\qquad w_c=p,\qquad z_c=g-r,\qquad y_c=y-t,
 \qquad s_c=s.
\end{equation}
The projection relations give $-z_c\in N_{\mathcal C}(p)$ and
$-y_c\in N_{\mathcal K}(s)$. They do not make a nonzero $F$ an optimality
certificate: the remaining discrepancies satisfy the exact identities
\begin{equation}\label{polish:eq:candidate-errors}
 \begin{aligned}
 Ap-b-s&=t-Ar,\\
 Qp+c-A^* y_c-z_c&=(I-Q)r+A^* t.
 \end{aligned}
\end{equation}
Operator magnitudes can amplify a small raw projected residual. Final native
feasibility, stationarity and objective gap must therefore be checked directly.

\subsection{Asynchronous completion}
Polishing starts from a promising HPR point, while the GPU iteration continues with its own anchor, penalty and restart controls. The correction matrix is held fixed during each polishing episode, allowing its linear-algebra work to be reused. An unsuccessful correction leaves the HPR solve running.

A candidate is accepted only after an independent check of primal feasibility, dual feasibility and relative objective gap in the original model, including finite-support and cone-domain conditions. The CPU correction does not feed intermediate iterates into HPR. It supplies an alternative final solution when it passes the requested accuracy test. Detailed triggering and linear-solve policies are documented with the implementation~\citep{hprqpc-release011}.

\section{Concluding remarks}
\label{sec:conclusions}
The common HPR framework separates the outer iteration from the design of its joint dual step. An sGS choice of $\mathcal T_1$ gives a sequential sweep; a block-linearized choice gives independent updates, with the second spectral moment tightening the row majorizer. Pairwise coupling extends the construction to convex QCQP while preserving simple cone projections. Equality Gram scaling improves row geometry, and asynchronous CPU polishing provides an independently checked route to high-accuracy completion.

The fixed-metric guarantees are established under the stated assumptions. Convergence under adaptive metric changes and the computational benefit of combining these techniques remain questions for further study.

\bibliographystyle{abbrv}
\bibliography{HPR-QP-arXiv/HPR-QP-ref,HPR-QP-arXiv/sections/polishing/references,HPR-QP-arXiv/sections/gram_scaling/references,HPR-QP-arXiv/sections/soc_support/references,HPR-QP-arXiv/sections/parallel/references}

\begin{thebibliography}{1}

\bibitem{parallel-hpr-framework}
K.~Chen, D.~Sun, Y.~Yuan, G.~Zhang, and X.~Zhao.
\newblock {HPR-QP}: A dual {Halpern Peaceman--Rachford} method for solving
  large-scale convex composite quadratic programming.
\newblock arXiv preprint, 2025.
\newblock \href{https://arxiv.org/abs/2507.02470v1}{arXiv:2507.02470v1}.
  Proposition, theorem and section numbers refer to this version.

\bibitem{scaling-curtis1972}
A.~R. Curtis and J.~K. Reid.
\newblock On the automatic scaling of matrices for {Gaussian} elimination.
\newblock {\em Journal of the Institute of Mathematics and Its Applications},
  10(1):118--124, 1972.
\newblock
  \href{https://doi.org/10.1093/imamat/10.1.118}{doi:10.1093/imamat/10.1.118}.

\bibitem{hprqpc-release011}
{HPR-QP-C contributors}.
\newblock {HPR-QP-C}, public release v0.1.1.
\newblock
  \href{https://github.com/PolyU-IOR/HPR-QP-C/tree/5f04c18dc45210f7d481901aac6bdc2f78414df4}{Source
  code and documentation}, 2026.

\bibitem{gram-coordinate-note}
{HPR-QP-C experimental branch}.
\newblock Equality {Gram} scaling in {HPR-QP}: Coordinate equivalence, {ADMM}
  interpretation, and implementation contract.
\newblock
  \href{https://github.com/PolyU-IOR/HPR-QP-C/blob/5f04c18dc45210f7d481901aac6bdc2f78414df4/docs/equality_gram_scaling.tex}{HPR-QP-C
  v0.1.1 documentation}, 2026.

\bibitem{li2018qsdpnal}
X.~Li, D.~Sun, and K.-C. Toh.
\newblock {Q}{S}{D}{P}{N}{A}{L}: {A} two-phase augmented {L}agrangian method
  for convex quadratic semidefinite programming.
\newblock {\em Math. Program. Comput.}, 10:703--743, 2018.
\newblock
  \href{https://doi.org/10.1007/s12532-018-0137-6}{doi:10.1007/s12532-018-0137-6}.

\bibitem{li2019block}
X.~Li, D.~Sun, and K.-C. Toh.
\newblock A block symmetric {G}auss--{S}eidel decomposition theorem for convex
  composite quadratic programming and its applications.
\newblock {\em Math. Program.}, 175:395--418, 2019.
\newblock
  \href{https://doi.org/10.1007/s10107-018-1247-7}{doi:10.1007/s10107-018-1247-7}.
  arXiv version: \href{https://arxiv.org/abs/1703.06629v2}{1703.06629v2}.

\bibitem{pock2011diagonal}
T.~Pock and A.~Chambolle.
\newblock Diagonal preconditioning for first order primal-dual algorithms in
  convex optimization.
\newblock In {\em Proceedings of the {IEEE} International Conference on
  Computer Vision}, pages 1762--1769. IEEE, 2011.
\newblock
  \href{https://doi.org/10.1109/ICCV.2011.6126441}{doi:10.1109/ICCV.2011.6126441}.

\bibitem{ruiz2001scaling}
D.~Ruiz.
\newblock A scaling algorithm to equilibrate both rows and columns norms in
  matrices.
\newblock Technical Report RAL-TR-2001-034, Rutherford Appleton Laboratory,
  2001.
\newblock
  \href{https://epubs.stfc.ac.uk/manifestation/271/raltr-2001034.pdf}{Public
  report}.

\bibitem{sun2025accelerating}
D.~Sun, Y.~Yuan, G.~Zhang, and X.~Zhao.
\newblock Accelerating preconditioned {ADMM} via degenerate proximal point
  mappings.
\newblock {\em SIAM J. Optim.}, 35(2):1165--1193, 2025.
\newblock \href{https://doi.org/10.1137/24M1650053}{doi:10.1137/24M1650053}.
  arXiv version: \href{https://arxiv.org/abs/2403.18618v2}{2403.18618v2}.

\end{thebibliography}
\end{document}